\documentclass[a4paper]{article}

\usepackage[margin=1in]{geometry}

\usepackage[T1]{fontenc}
\newcommand{\changefont}[3]{
\fontfamily{#1} \fontseries{#2} \fontshape{#3} \selectfont}

\changefont{ptm}{m}{n}

\usepackage{setspace} 
\usepackage{graphicx}    

\usepackage[all]{xy,xypic}
\usepackage{amsfonts,amssymb,amsmath,amsgen,amsopn,amsbsy,theorem,graphicx,epsfig}
\usepackage{eufrak,amscd,bezier,latexsym,mathrsfs,enumerate}\usepackage[utf8]{inputenc}
\usepackage[english]{babel}
\usepackage[dvipsnames]{xcolor}

\usepackage{amssymb}
\usepackage{graphics}
\usepackage{mathrsfs}
\usepackage{color}

\newtheorem{theorem}{Theorem}[section]

\newtheorem{lemma}{Lemma}[section]

\newtheorem{definition}{Definition}[section]

\long\def\symbolfootnote[#1]#2{\begingroup%
\def\thefootnote{\fnsymbol{footnote}}\footnote[#1]{#2}\endgroup} 

\begin{document}

\begin{center}
\Large \textbf{Unpredictable Solutions of Linear Differential Equations}
\end{center}

\begin{center}
\normalsize \textbf{Marat Akhmet$^{1,}\symbolfootnote[1]{Corresponding Author Tel.: +90 312 210 5355,  Fax: +90 312 210 2972, E-mail: marat@metu.edu.tr}$, Mehmet Onur Fen$^2$, Madina Tleubergenova$^{3,4}$, Akylbek Zhamanshin$^{3,4}$} \\
\vspace{0.2cm}
\textit{\textbf{\footnotesize$^1$Department of Mathematics, Middle East Technical University, 06800, Ankara, Turkey}} \\
\textit{\textbf{\footnotesize$^{2}$Department of Mathematics, TED University, 06420, Ankara, Turkey}} \\
\textit{\textbf{\footnotesize$^{3}$Department of Mathematics, Aktobe Regional State University, Aktobe, Kazakhstan}} \\
\textit{\textbf{\footnotesize$^{4}$Institute of Information and Computational Technologies, Almaty, Kazakhstan}}
\vspace{0.1cm}
\end{center}

\vspace{0.3cm}

\begin{center}
\textbf{Abstract}
\end{center}

\noindent\ignorespaces

In this study, the existence and uniqueness of the unpredictable solution for a non-homogeneous linear system of ordinary differential equations is considered. The hyperbolic case is under discussion. New properties of unpredictable functions are discovered. The presence of the solutions confirms the existence of Poincar\'e chaos. Simulations illustrating the chaos are provided.

\vspace{0.2cm}
 
\noindent\ignorespaces \textbf{Keywords:} Unpredictable solutions; Poincar\'e chaos; Linear non-homogeneous systems

\vspace{0.6cm}


\section{Introduction and Preliminaries}\label{neural_prelim}

The concept of the unpredictable functions was introduced in the paper \cite{a1}. The description of such functions relies on the dynamics of unpredictable points, which were presented in the study \cite{a2} for the first time in the literature. Considering unpredictable points the authors extended the limits of the theory of classical dynamical systems, which was founded by H. Poincar\'e and G. Birkhoff. The definition in \cite{a1} considers the functions as points of the Bebutov dynamics \cite{sell}. The metric of the dynamics is not convenient for  applications in the theory of differential equations. Therefore, in paper \cite{a4}, it was suggested to utilize the topology of uniform convergence on compact subsets of the real axis. This definition is possibly the most effective one for methods of qualitative analysis. The paper \cite{a4} was devoted to the investigation of sufficient conditions for the existence of unpredictable solutions of quasilinear differential equations in the case that matrices of coefficients admit all eigenvalues with negative real parts as well as discrete equations. 

The present study has two principal novelties with respect to the previous results in the field. The first one is that we consider the hyperbolic case, when the eigenvalues of the matrix of coefficients can admit positive real parts. The second one is that we propose a simpler and more comprehensible proof for the unpredictability property this time. As it was confirmed in papers by Akhmet and Fen \cite{a2}-\cite{a4}, the existence of unpredictable solutions simultaneously means the presence of Poincar\'e chaos, i.e., unpredictable solutions are ``irregular".  This makes the subject attractive for applications. Finally, we consider new properties of the functions. Sufficient conditions are provided such that a linear transformation of an unpredictable function is unpredictable, and it is proved that the sum of an unpredictable function and a periodic function is an unpredictable function.

In the remaining part of the paper, we will make use of the usual Euclidean norm for vectors and the norm induced by the Euclidean norm for square matrices.

The definition of unpredictable functions is as follows. 

\begin{definition}  \label{def1}  \cite{a4} A uniformly continuous and bounded function $\vartheta: \mathbb R\rightarrow \mathbb R^m$ is unpredictable if there exist positive numbers $\epsilon_0$, $\delta$ and sequences $\{t_n\}$, $\{u_n\}$ both of which diverge to infinity such that  $\|\vartheta(t+t_n)-\vartheta(t)\|\rightarrow0$ as $n\rightarrow\infty$ uniformly on compact subsets of the axis and $\|\vartheta(t+t_n)-\vartheta(t)\|\geq\epsilon_0$ for each $t\in[u_n-\delta, u_n+\delta]$ and $n\in\mathbb{N}$.
\end{definition}

For the convenience of the next  discussion, we will call  the convergence of the  function's shifts on compact  subsets  as \textit{Poisson stability} and  the existence of  the number $\epsilon_0$ as \textit{unpredictability property} of the function.
Thus,  a function is \textit{unpredictable}, if it  is \textit{Poisson stable} and admits  the \textit{unpredictability property}.   

It is worth noting that in the literature a large number of results are obtained for periodic, quasi-periodic and almost periodic solutions of differential equations due to the established mathematical methods and important applications. On the other hand, recurrent and Poisson stable solutions are also crucial for the theory of differential equations \cite{bender,sell}. The proposal by Akhmet and Fen can revive interest of specialists in differential equations theory for two reasons. The first one is related to the verification of the unpredictability which  requests more sophisticated technique than for recurrent and Poisson stable solutions. Thus the problem of the existence of unpredictable solutions is a challenging one. Another crucial reason is that the presence of unpredictable solutions is necessarily accompanied by the existence of chaos, which is called \textit{Poincar\'e chaos} in papers \cite{a2}-\cite{a4}. Consequently, the research of differential equations with unpredictable solutions will definitely help to activate study and applications of chaos.    

The main object of the present paper is the following system of linear differential equations,
\begin{equation}\label{integral_eqn}
x'(t)=  Ax(t)+g(t),
\end{equation} 
where $x \in \mathbb R^n,$ $n$ is a fixed natural number and $\mathbb R$ is the set of all real numbers. Moreover, we assume that all eigenvalues of the constant matrix $A\in \mathbb R^{n \times n}$ have nonzero real parts, and the function $g:\mathbb R \to \mathbb R^n$ is uniformly continuous and bounded. 

Assume that  $\Re e \lambda_i  < 0,  i =1,2,\ldots,q,$ and $\Re e \lambda_i > 0, i = q+1,\ldots,p,  1 \le q <p,$ where  $\lambda_i, i = 1,\ldots,p,$  are  eigenvalues of the matrix $A.$  It  is known that  one can  find a regular  matrix $B$  such  that  the transformation $x = By$   changes the system (\ref{integral_eqn})   to  the equation 

\begin{equation}\label{integral_eqn1}
y'(t)=  B^{-1}ABy(t)+B^{-1}g(t),
\end{equation} 
with the diagonal  matrix of coefficients.   This is why, without  loss of generality one can assume that the matrix $A$ is equal to $\mbox{diag}(A_-,A_+),$ where eigenvalues of  matrices $A_-$ and $A_+$  are  with  negative  and positive real  parts respectively.  From the equation (\ref{integral_eqn1}), it  implies that  the following auxiliary   assertion will be needed.
\begin{lemma} \label{lem1} The function $f(t)= B^{-1}g(t)$ is unpredictable.
\end{lemma} 
The proof of the lemma immediately follows the  inequalities  $\| f(t+t_n) - f(t)\|  \le \|B^{-1}\|\|	g(t+t_n) - g(t)\|$ and $\| f(t+t_n) - f(t)\| \ge \frac{1}{\|B\|}\|	g(t+t_n) - g(t)\| .$  

For further applications of the main result of the paper, the following lemma can be useful.

\begin{lemma} \label{lem2} Assume that $g(t)$ is an unpredictable function, and a function $f(t)$ is continuous and periodic. Then the sum $g+f$ is an unpredictable function. 
\end{lemma}	 
 \textbf{Proof.}  Consider the sum $h(t)\equiv g(t) + f(t).$ Let  $t_n$  be the sequence such that $g(t+t_n)$ converges to  $g(t)$  uniformly  on all  compact  subsets  of the axis.  Since of periodicity  of $f(t),$ one can find a subsequence of $t_n$ (we will assume, without  loss of generality that it  is the sequence $t_n$ itself) that satisfy $f(t+t_n) - f(t) \to 0$ uniformly on the axis. Consequently,  $h(t+t_n)-h(t) = (g(t+t_n) - g(t)) + (f(t+t_n)-f(t))  \to  0$  uniformly  on compact  subsets of the axis  as $n \to \infty.$ The Poisson stability  is proved.
 
Consider, now, the sequence $u_n \to  \infty$ and positive numbers $\eta, \sigma,$ such that $\|g(t+t_n) - g(t)\| \ge \eta$  for $t \in[u_n-\sigma,u_n + \sigma].$   Again, due to  the periodicity,  one can  find a subsequence of $u_n,$  let  say  the sequence itself,  such  that $\|f(t+t_n)-f(t)\| < \eta/2$ if $t \in[u_n-\sigma,u_n + \sigma].$   Then,  $\|h(t+t_n)-h(t)\| \ge \|g(t+t_n) - g(t)\| - \|f(t+t_n)-f(t)\| > \eta -\eta/2 = \eta/2$  for  $t \in[u_n-\sigma,u_n + \sigma].$  The unpredictability  property  is proved. $\Box$ 

\section{Main result}\label{main}

In what follows we will denote $g=(g_-,g_+),$ where the vector-functions $g_-$ and $g_+$ are of dimensions $q$ and $p-q$, respectively.
    
As it is known from the theory of differential equations, system (\ref{integral_eqn}) admits a unique solution $\varphi(t) = (\varphi_-(t),\varphi_+(t))$ which is bounded on $\mathbb R$, where 
\begin{equation}\label{dif_eqn}
\varphi_-(t)=\int_{-\infty}^t e^{A_-(t-s)}g_-(s)ds, \quad   \varphi_+(t)=-\int_{t}^{\infty} e^{A_+(t-s)}g_+(s)ds.
\end{equation} 
One can confirm that the bounded solution is periodic, quasi-periodic, or almost periodic if the perturbation function $g$ is respectively of the same type. 

The following theorem is concerned with the unpredictable solution of system  (\ref{integral_eqn}).
\begin{theorem} \label{mainthm} If all   eigenvalues of the matrix $A$ admit non-zero real  parts, and the   function $g(t)$ is unpredictable, then the  system   (\ref{integral_eqn})  possesses a unique  unpredictable solution. If additionally,  all the eigenvalues are  with  negative real parts, then the unpredictable solution  is uniformly  asymptotically  stable. 
\end{theorem}
\textbf{Proof.}  From the conditions of the theorem  it  implies that  the system  admits a bounded solution and it  is uniformly  asymptotically  stable if the matrix of coefficients admits all  eigenvalues with negative real parts.  This is why, it is sufficient to prove that the solution (\ref{dif_eqn}) is an unpredictable function.

We   will show,  first,  that the solution  is Poisson stable.   Since the eigenvalues of the matrix A in system  (\ref{integral_eqn})  have   non-zero  real parts, there exist numbers $ K\geq1$ and $\alpha>0$ such that $\|e^{A_-t}\|\leq Ke^{-\alpha t}$ for $t\geq0$  and  $\|e^{A_+t}\|\leq Ke^{\alpha t}$ for $t\leq0.$ 

Since   the  function $g$   is Poisson stable there  exists a sequence $t_n \to \infty$   such  that  $\|g(t+t_n) - g(t)\| \to  0$   uniformly  on compact  subsets  of the  axis.  
One can  easily  find that 
\begin{equation*}\begin{split}
    & \|\varphi_-(t+t_n)-\varphi_-(t)\|= \|\int_{-\infty}^t e^{A_-(t-s)}g_-(s+t_n)ds- 
    \int_{-\infty}^t e^{A_-(t-s)}g_-(s)ds\| \\
	& = \|\int_{-\infty}^t e^{A_-(t-s)}[g_-(s+t_n)-g_-(s)]ds\| \leq  \int_{-\infty}^t\|e^{A_-(t-s)}\|\|g_-(s+t_n) \\
	& -g_-(s)\|ds  \leq \int_{-\infty}^t Ke^{-\alpha(t-s)}\|g_-(s+t_n)-g_-(t)\|ds,
	\end{split}
\end{equation*}
and 
\begin{equation*}\begin{split}
	& \|\varphi_+(t+t_n)-\varphi_+(t)\|=\|\int_{t}^{\infty} e^{A_+(t-s)}g_+(s+t_n)ds-\int_{t}^{\infty} e^{A_+(t-s)}g_+(s)ds\| \\
	& =\|\int_{t}^{\infty}e^{A_+(t-s)}[g_+(s+t_n)-g_+(s)]ds\|\leq\int_{t}^{\infty}\|e^{A_+(t-s)}\|\|g_+(s+t_n) \\
	& -g_+(s)\|ds\leq	\int_{t}^{\infty}Ke^{\alpha(t-s)}\|g_+(s+t_n)-g_+(t)\|ds.
\end{split}\end{equation*}

Fix an  arbitrary  positive number $\epsilon$  and a section $[a,b],  -\infty <a<b<\infty,$  of the  real  axis.
We will show that  for  sufficiently  large $n$  it  is true that  $\|\varphi(t+t_n)-\varphi(t)\| <\epsilon$  on $[a,b].$   Denote  $M =\sup_{\mathbb R}\|g\|,$  and choose numbers $c <a, b <d,  \xi > 0,$  such  that  
$\frac{2MK}{\alpha}e^{-\alpha(a-c)} < \frac{\epsilon}{4},   \frac{2MK}{\alpha}e^{-\alpha(d-b)} < \frac{\epsilon}{4}  $ and  $ \frac{K\xi}{\alpha} <  \frac{\epsilon}{4}.$

Consider $n$   sufficiently  large such  that  $\|g(t+t_n) - g(t)\| < \xi$  on $[c,b].$ Then for  all  $t \in [a,b]$  
\begin{equation*}\begin{split}
	&\|\varphi_-(t+t_n)-\varphi_-(t)\|\leq \int_{-\infty}^{c}K e^{-\alpha(t-s)}\|g_-(s+t_n)-g_-(s)\|ds+ \\
	&\int_{c}^t Ke^{-\alpha(t-s)}\|g_-(s+t_n)-g_-(s)\|ds\leq\int_{-\infty}^{c}K e^{-\alpha(t-s)}2Mds+ \\
	&\int_{c}^t Ke^{-\alpha(t-s)}\xi ds\leq 
	\frac{2MK}{\alpha}e^{-\alpha(a-c)} + \frac{K\xi}{\alpha} \leq  \frac{\epsilon}{4} + \frac{\epsilon}{4} = \frac{\epsilon}{2},
\end{split}\end{equation*}
and similarly one can  show that    
\begin{equation*}\begin{split}
	&\|\varphi_+(t+t_n)-\varphi_+(t)\|\leq \int_{t}^{d}K e^{\alpha(t-s)}\|g_+(s+t_n)-g_+(s)\|ds+ \\
	&\int_{d}^{\infty} Ke^{\alpha(t-s)}\|g_+(s+t_n)-g_+(s)\|ds\leq\int_{t}^{d}K e^{\alpha(t-s)}\xi ds+ \\
	&\int_{d}^t Ke^{-\alpha(t-s)} 2M ds\leq 
	\frac{K\xi}{\alpha}  + \frac{2MK}{\alpha}e^{-\alpha(d-b)} \le   \frac{\epsilon}{4} + \frac{\epsilon}{4} = \frac{\epsilon}{2}.
\end{split}\end{equation*}

We have obtained that  for  sufficiently  large  $n$  it  is true that  $$\|\varphi(t+t_n)-\varphi(t)\|  
\le  \|\varphi_+(t+t_n)-\varphi_+(t)\| + \|\varphi_-(t+t_n)-\varphi_-(t)\| < \epsilon$$
for $t \in [a,b].$ 

The Poisson stability  of the  solution $\varphi$  is proved. 

Let  us check  that  the  solution possesses the  unpredictability  property.   Since the  function $g$  is unpredictable,  there  exist  a  sequence $u_n  \to  \infty$ and  positive numbers $\eta,$ $\kappa$ such  that $\|g(u_n+t) - g(u_n)\| \geq \eta$  for 
 $t \in [u_n - \kappa, u_n+\kappa].$    One   can  easily  check that either $\|g_-(u_n+t) - g_-(u_n)\| \geq \eta/2$  or $\|g_+(u_n+t) - g_+(u_n)\| \geq \eta/2$  for  the same $t.$   Assume that  the first  inequality  is valid, since another case can    be considered very  similarly.  Then, there  is a number, without  loss of generality  $\kappa$ itself,  such  that  $\|g_-(t+t_n) - g_-(t)\| \geq \eta/4$ for  all  $t \in [u_n - \kappa, u_n+\kappa].$  Without  loss of generality  we assume that  $\|g_-(t+t_n) - g_-(t)\| \geq \eta$ for  all $t \in [u_n - \kappa, u_n+\kappa].$

At  first we assume that   $\inf_n\|\varphi(u_n+t_n)-\varphi(u_n)\|\ =0.$  This contradicts the equality  
$$\varphi(t+t_n)-\varphi(t) =   \varphi(u_n+t_n)-\varphi(u_n) + $$$$ \int_{u_n}^t  A(\varphi(s+t_n)-\varphi(s)) ds  
+   \int_{u_n}^t (g(s+t_n)-g(s))ds,$$
since    the last  integral  is a positive number larger than $t\eta$  and other terms tend to  zero   as $n \to  \infty$   for   each  positive fixed moment  $t$  from the interval  $ [u_n - \kappa, u_n+\kappa].$

This is why, it  is true that     $\inf_n\|\varphi(u_n+t_n)-\varphi(u_n)\|\ =2\epsilon_0,$  where  $\epsilon_0$ is a positive number.   Fix a positive number  $\sigma$  such  that  $2\eta\sigma e^{\|A\|\kappa} < \epsilon_0.$ 

Now, we have that  
$$\|\varphi(t+t_n)-\varphi(t)\|= \|\varphi(u_n+t_n)-\varphi(u_n) + $$$$ \int_{u_n}^t e^{A(t-s)}[g(s+t_n)-g(s)]ds\|\geq 2\epsilon_0 -  \eta\sigma e^{\|A\|\kappa} > \epsilon_0, $$
for   $t  \in [u_n -\sigma, u_n+\sigma].$ 
The  unpredictability  is proved.  

The proof of the theorem is completed. $\Box$

\section{Examples}\label{Example}    

One of the possible ways to confirm the presence of chaos is through simulations. The concept of unpredictable solutions maintain the series of oscillators, but from the other side the chaos accompanies unpredictability. Consequently, we can look for a confirmation of the results   for unpredictability observing irregularity in simulations. The approach is effective for asymptotically stable unpredictable solutions, and it is just illustrative for hyperbolic systems with unstable solutions. In the latter case we rely on the fact that any solution becomes unpredictable ultimately.    

In the following examples we will utilize the function
 \begin{eqnarray}  \label{func_unprdctble_psi}
 \Theta(t)=\displaystyle \int_{-\infty}^{t} e^{-2(t-s)} \Omega(s) ds,
 \end{eqnarray} 
 which was discussed in paper \cite{a4}.   The function  $\Omega(t)$ is defined by $\Omega(t)=\psi_{i}$ for $t\in [i,i+1),$ $i\in\mathbb Z, $  where  $\{\psi_{i}\},$ $i\in \mathbb Z,$ is an unpredictable solution of the logistic map
 \begin{eqnarray} \label{logistic_map_example}
 \lambda_{i+1}= 3.91 \lambda_{i}(1-\lambda_{i})
 \end{eqnarray}
 inside the unit interval $[0,1].$  
 The function $\Theta(t)$ is bounded on the whole real axis such that $\displaystyle \sup_{t\in\mathbb R} \left|\Theta(t)\right| \leq 1/2$, and it is uniformly continuous since its derivative is bounded.
 It was proven in paper \cite{a4} that $\Theta(t)$ is an unpredictable function.\\ 
 
\noindent \textbf{Example 3.1.} 

\begin{figure}[ht]
	\centering
	\includegraphics[height=7.5cm]{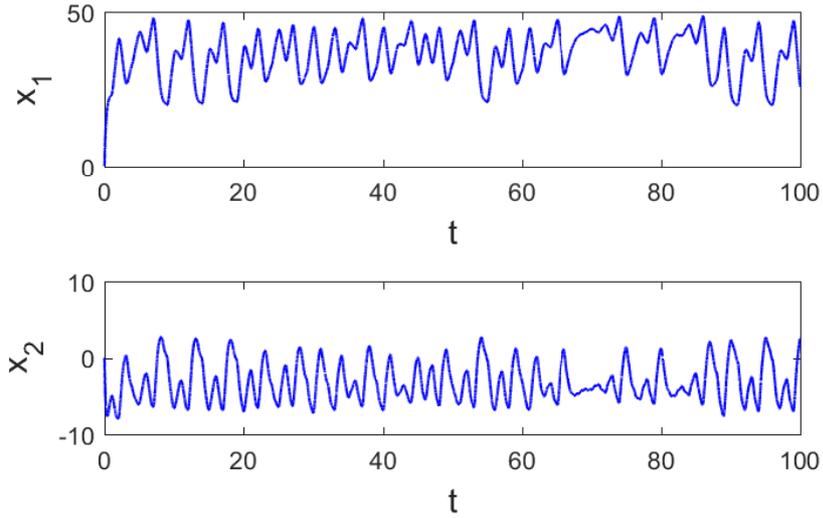}
	\caption{The time series of the $x_1$ and $x_2$ coordinates of system (\ref{2}) with the initial conditions $x_1(0)=0,18,$ $x_2(0)=0,01.$ The figure manifests the presence of Poincar\'e chaos.}
	\label{2.5.4}
\end{figure}

\begin{figure}[ht]
	\centering
	\includegraphics[height=6.0cm]{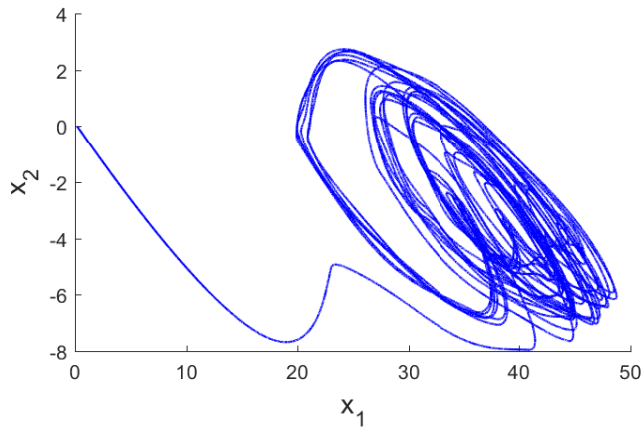}
	\caption{The trajectory of system (\ref{2}) with the initial point $(0,0,01).$}
	\label{2.5.5}
\end{figure}

Consider the system 
\begin{eqnarray} 
\label{2}
\begin{array}{l}
x_1'=-2x_1+2x_2+259\Theta(t)-\sin(10t)\\
x_2'=x_1-3x_2-150\Theta(t)+\cos(10t)
\end{array}
\end{eqnarray}
where the eigenvalues of the matrix of coefficients are $-2$  and $-0.5$.  One can confirm that the perturbation function $(259\Theta(t)-\sin(10t), -150\Theta(t)+\cos(10t))$ is unpredictable in accordance with Lemma \ref{lem2}. By the main result of our paper, there  is an  asymptotically stable unpredictable solution $(\varphi_1(t),\varphi_2(t))$ of system (\ref{2}). Consequently, any solution of the equation behaves irregularly ultimately. This is seen from the simulation of the solution with $x_1(0)=0,18$, $x_2(0)=0,01$ in Figures \ref{2.5.4} and \ref{2.5.5}.

The next example is devoted to a system with hyperbolic linear part such that the matrix of coefficients admit both positive and negative eigenvalues. 

\newpage

\noindent \textbf{Example 3.2.}  

Let us take into account the the system
\begin{eqnarray} \label{3b}
\begin{array}{l}
u_1'=-52098 u_1+7090 \varphi_2(t) \\
u_2'=9.5 u_1+0.0000000325 u_2+0.111 \varphi_1(t),
\end{array}
\end{eqnarray}
where $(\varphi_1(t),\varphi_2(t))$ is the unpredictable solution of system (\ref{2}). 

The eigenvalues of the matrix of coefficients  of system (\ref{3b}) are $-52098$ and $0.0000000325.$ According to the result of Theorem \ref{mainthm}, system (\ref{3b}) possesses a unique unpredictable solution.

\begin{figure}[ht]
	\centering
	\includegraphics[height=7.5cm]{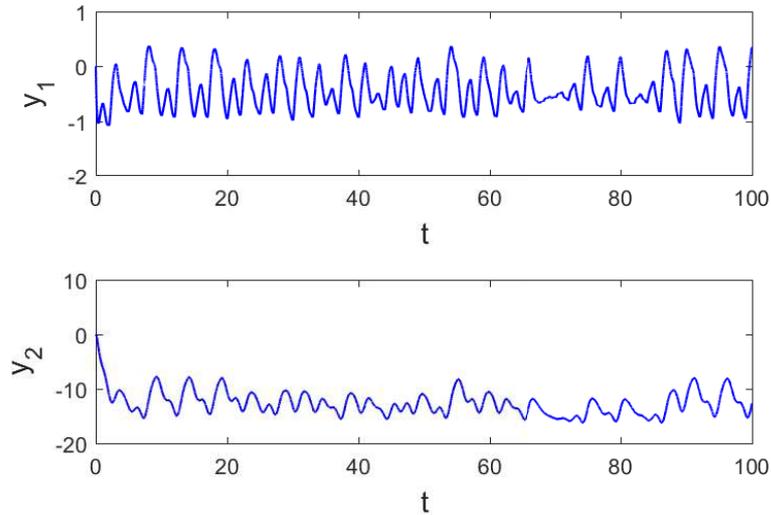}
	\caption{The time series for the $y_1$ and $y_2$ coordinates of system (\ref{3}) with the initial conditions $y_1(0)=0,$ $y_2(0)=0.$ The irregular behavior of the solution reveals the presence of an unpredictable solution in the dynamics of (\ref{3}).}
	\label{2.5.6}
\end{figure}

\begin{figure}[ht]
	\centering
	\includegraphics[height=6.0cm]{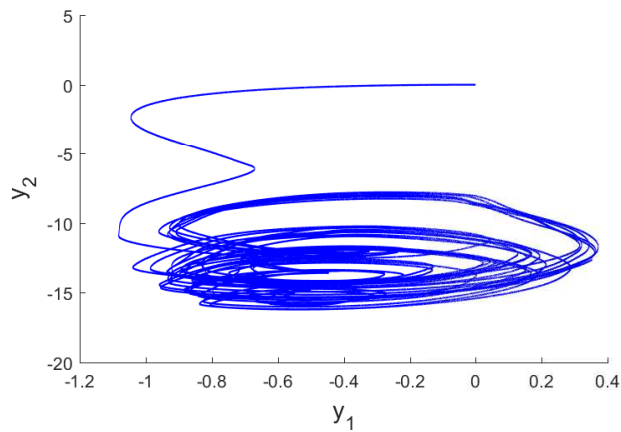}
	\caption{The trajectory of the solution of system (\ref{3}) with $y_1(0)=0$ and $y_2(0)=0.$}
	\label{2.5.7}
\end{figure}

In order to show the irregular behavior, we consider the system
\begin{eqnarray} 
\label{3}
\begin{array}{l}
y_1'=-52098y_1+7090x_2(t)\\
y_2'=9.5y_1+0.0000000325y_2+0.111x_1(t),
\end{array}
\end{eqnarray}
where $(x_1(t),x_2(t))$ is the solution of (\ref{2}) depicted in Figures \ref{2.5.4} and \ref{2.5.5}. The simulation results for system (\ref{3}) corresponding to the initial conditions  $y_1(0)=0$ and $y_2(0)=0$ are shown in Figures \ref{2.5.6} and \ref{2.5.7}. Both of the figures confirm the presence  of  unpredictability in the dynamics of system (\ref{3}).

\section{Acknowledgement}

The third and fourth authors were supported in parts by the MES RK grant
No. AP05132573 ``Cellular neural networks with continuous/discrete time and singular
perturbations." (2018-2020)
of the Committee of Science, Ministry of Education and Science of the
Republic  of Kazakhstan.

\end{document}